\documentclass[reqno,10pt,a4paper]{amsart}
\usepackage{srcltx}
\usepackage{graphicx,color,amssymb,latexsym,amsmath,amsfonts}
\usepackage{mathrsfs}
\numberwithin{equation}{section}
\theoremstyle{plain}
\newtheorem{keyword}{keyword}
\newtheorem{thm}{Theorem}[section]

\newtheorem{lemma}{Lemma}[section]

\newtheorem{proposition}{Proposition}[section]

\begin{document}

\title{On asymptotically minimax nonparametric detection\\
of signal in Gaussian white noise}

% indicate corresponding author with \corref{}
% \author{\fnms{John} \snm{Smith}\thanksref{t2}\corref{}\ead[label=e1]{smith@foo.com}\ead[label=e2,url]{www.foo.com}}
% \thankstext{t2}{Thanks to somebody}
% \address{line 1\\ line 2\\ \printead{e1}\\ \printead{e2}}
\
\author{Mikhail Ermakov}\\

\maketitle

{erm2512@gmail.com}\\
{Institute of Problems of Mechanical Engineering, RAS and\\
St. Petersburg State University, St. Petersburg, RUSSIA\\
{Mechanical Engineering Problems Institute\\
Russian Academy of Sciences\\
Bolshoy pr.,V.O., 61\\
St.Petersburg\\
Russia}}
\vskip 0.3cm
St.Petersburg State University\\
University pr., 28, Petrodvoretz\\
198504 St.Petersburg\\
Russia
\vskip 0.3cm
St.Petersburg department of\\ Steklov
mathematical institute \\
Fontanka 27, St.Petersburg 191023

\begin{abstract} For the problem of nonparametric detection of signal in Gaussian white noise we point out strong aasymptotically minimax tests. The sets of alternatives  are a ball   in Besov space $B^r_{2\infty}$ with "small" balls in $L_2$ removed.
\end{abstract}

\begin{keyword}[class=AMS]
[Primary ]{62G10},{62G20},{62M02}
%\kwd{}
%\kwd[; secondary ]{}
\end{keyword}

\begin{keyword}
 signal detection, asymptotic minimaxity, asymptotic efficiency, Gaussian white noise
%\kwd{}
\end{keyword}

% history:

%\tableofcontents

\section{\bf  Introduction. Main Result. }
Let we observe a random process $Y_\epsilon(t),t \in (0,1), \epsilon > 0$, defined by stochastic differential equation
\begin{equation}\label{i1}
dY_\epsilon(t) = \theta(t)\,dt + \epsilon dw(t)
\end{equation}
with Gaussian white noise $w(t)$. The signal $\theta \in L_2(0,1)$ is unknown.

Our goal is to  test the hypothesis
$$ H_0: \, \theta(t) = 0,\,\,\, t \in (0,1) $$
versus the alternative
$$H_\epsilon:\quad  \int_0^1 \theta^2(t)\,dt = ||\theta||^2 > \rho_\epsilon> 0,
$$
if a priori information is provided that
$$
\theta \in B^r_{2\infty}(P_0) =\left\{\theta:\, \theta(t)= \sum_{j=1}^\infty\theta_j\phi_j(t),\, k^{-2r}\sum_{j=k}^\infty\theta_j^2 \le P_0, t \in (0,1), 1 \le k < \infty\right\}$$
with $P_0 >  0$. Here $\phi_j, 1 \le j<\infty$, is orthonormal system of functions. For wide class of orthonormal systems of functions $\phi_j, 1 \le j<\infty$ the space
$$
\left\{\theta:\, \theta(t)= \sum_{j=1}^\infty\theta_j\phi_j(t),\, k^{-2r}\sum_{j=k}^\infty\theta_j^2  < \infty, t \in (0,1), 1 \le k < \infty\right\}$$ is Besov space $B^r_{2\infty}$ (see \cite{rio})

Denote $Q_\epsilon= \{\theta:\, ||\theta||^2 \ge \rho_\epsilon, \theta \in  B_{2\infty}^r(P_0)\}$.

For any test $K\epsilon$ denote $\alpha(K_\epsilon)$ its type I error probability and denote $\beta_\theta(K_\epsilon)$ its type II error probability  for the alternative $\theta \in Q_\epsilon$.

We put
$$
\beta_\epsilon(K_\epsilon) = \sup_{\theta\in Q_\epsilon} \beta_\theta(K_\epsilon).
$$
We say that family of tests $L_\epsilon$ is asymptotically minimax if, for any family of tests $K_\epsilon$, $\alpha(K_\epsilon) \le \alpha(L_\epsilon)$, there holds
$$
\lim\sup_{\epsilon \to 0} \beta_\epsilon(K_\epsilon) - \beta_\epsilon(L_\epsilon) \ge 0.$$
Paper goal is to establish asymptotically minimax families of tests $L_\epsilon$ for the sets of alternatives $Q_\epsilon$. If the sets of alternatives are ellipsoids with "small balls" removed,
asymptotically minimax families of tests have been found in \cite{er90}. For nonparametric hypothesis testing this result can be considered as a version  of Pinsker Theorem \cite{pin,ts,jo} on asymptotically minimax nonparametric estimation. Note that hypothesis testing with nonparametric sets of alternatives belonging some ball in functional space is intensively studied (see \cite{is, dal} and references therein).

The proof, in main features, repeats the reasoning in \cite{er90}. The main difference  in the proof is the solution of another extremal problem minimizing type II error probabilities
caused another definition of sets of alternatives. Other differences have technical character and are also caused the differences of definitions of sets of alternatives.

Define $k=k_\epsilon$ and $\kappa^2 = \kappa_\epsilon^2$ as a solution of two equations
\begin{equation}\label{i2}
2rk_\epsilon^{2r+1}\kappa^2_\epsilon = P_0
\end{equation}
and
\begin{equation}\label{i3}
k_\epsilon\kappa_\epsilon^2 + k_\epsilon^{-2r} P_0  = \rho_\epsilon.
\end{equation}
Denote
$
\kappa_j^2 =  \kappa_\epsilon^2$, for $1 \le j \le k_\epsilon
$ and $\kappa_j^2 P_0(2r)^{-1} j^{-2r-1}$, for $ j > k_\epsilon.$

Define test statistics
$$
T^a_\epsilon(Y_\epsilon) =   \epsilon^{-4} \sum_{j=1}^\infty \kappa_j^2 y_j^2.
$$e
$$
A_\epsilon = \epsilon^{-4}\sum_{j=1}^\infty \kappa_j^4.
$$
For type I error probabilities $\alpha, 0 < \alpha< 1,$ define critical regions
$$
S^a_\epsilon == \{ y: \, (T^a_\epsilon(y) - \epsilon^{-2}\rho_\epsilon)(2A_\epsilon)^{-1/2} > x_\alpha\}
$$
with $x_\alpha$ defined by equation $$\alpha = 1 - \Phi(x_\alpha) = (2\pi)^{-1/2}\int_{x_\alpha}^\infty  \exp\{-t^2/2\} \, dt.$$
\begin{thm}\label{t1} Let
\begin{equation}\label{i4}
0 <\lim\inf_{\epsilon\to 0} A_\epsilon \le <\lim\sup_{\epsilon\to 0} A_\epsilon< \infty.
\end{equation}
Then the tests $L^a_\epsilon$ with critical regions $S^a_\epsilon$ are asymptotically minimax with $\alpha(L^a_\epsilon) = \alpha(1+o(1))$ and
 \begin{equation}\label{i5}
\beta_\epsilon(L_\epsilon^a) = \Phi(x_\alpha - (A_\epsilon/2)^{1/2})(1 + o(1))
\end{equation}
as $\epsilon \to 0$.
\end{thm}
\noindent{\sl Example.} Let $\rho_\epsilon = R\epsilon^{\frac{8\beta}{4\beta+1}}.$ Then
$$
A_\epsilon = \left(\frac{P_0}{2r}\right)^{1/2r}\frac{4r+2}{4r+1}\left(\frac{R}{2r+1}\right)^{\frac{4r-1}{2r}}.
$$
In what follows, we shall denote letter $C$ and $C$ with indices different generic constants.
\section{\bf Proof of Theorem \ref{t1}}
Fix $\delta, 0<\delta<1.$ Denote $\kappa_j^2(\delta) = 0$ for $j > \delta^{-1}k_\epsilon$. Define $\kappa_j^2(\delta), 1 \le j < k_{\epsilon\delta}=\delta^{-1}k_\epsilon,$  the equations (\ref{i2}) and (\ref{i3}) with $P_0$ and $\rho_\epsilon$ replaced with $P_0(1-\delta)$ and $\rho_\epsilon(1 + \delta)$ respectively.
Similarly to \cite{er90}, we find Bayes test for a priori distribution $\theta_j = \eta_j=\eta_j(\delta), 1 \le j < \infty,$ with Gaussian independent random variables $\eta_j, E\eta_j =0, E\eta_j^2 = \kappa_j^2(\delta)$ and show that this test is asymptotically minimax for some $\delta=\delta_\epsilon \to 0$ as $ \epsilon \to 0$.
\begin{lemma}\label{l3} For any $\delta,  0 < \delta <1,$  there holds
\begin{equation}\label{i7}
P(\eta(\delta) = \{\eta_j(\delta)\}_{j=1}^\infty  \in Q_\epsilon) = 1 + o(1)
\end{equation}
as $\epsilon \to 0$.
\end{lemma}
Denote
$$  A_{\epsilon,\delta} = \epsilon^{-4}\sum_{j=1}^\infty \kappa_j^4(\delta).
$$
By straightforward calculations, we get
\begin{equation}\label{i8}
\lim_{\delta \to 0}\lim_{\epsilon \to 0} A_\epsilon A_\epsilon^{-1}(\delta) =1.
\end{equation}
Denote $\gamma_j^2(\delta) = \kappa_j^2(\delta)(\epsilon^2 + \kappa_j^2(\delta))^{-1}$.

By Neymann-Pearson Lemma the Bayes critical region is defined the inequality
\begin{equation}\label{i9}
\begin{split}&
C_1 < \prod_{j=1}^{k_{\epsilon\delta}}(2\pi)^{-1/2}\kappa_j^{-1}(\delta) \int \exp\left\{- \sum_{j=1}^{k_{\epsilon\delta}}(2\gamma_j^2(\delta))^{-1}(u_j- \gamma_j^2(\delta)y_j)^2\right\}\exp\{-T_{\epsilon\delta}(y)\}\\&
= \exp\{-T_{\epsilon\delta}(y)\}(1+o(1))
\end{split}
\end{equation}
where
$$
T_{\epsilon\delta}(y) = \epsilon^{-2}\sum_{j=1}^\infty \gamma_j^2(\delta)y_j^2.
$$ 
Define critical region
$$
S_{\epsilon\delta}=\{y :  \, R_{\epsilon\delta}(y) = (T_{\epsilon\delta}(y) -C_{\epsilon\delta})(2A_\epsilon(\delta))^{-1/2} > x_\alpha\}
$$
with
$$
C_{\epsilon\delta} = E_0 T_{\epsilon\delta}(y) = \epsilon^{-2}\sum_{j=1}^\infty \gamma_j^2(\delta).
$$
Denote $L_{\epsilon\delta}$ the tests with critical regions $S_{\epsilon\delta}$.

Denote $\gamma_j^2 = \kappa_j^2(\epsilon^2 + \kappa_j^2)^{-1}, 1 \le j < \infty$
 Define test statistics $T_{\epsilon}, R_\epsilon$,  critical regions $S_\epsilon$ and constants $C_\epsilon$  by the same way as test statistics $T_{\epsilon\delta}, R_{\epsilon\delta}$, critical regions $S_{\epsilon\delta}$ and constants $C_{\epsilon,\delta}$ respectively with $\gamma_j^2(\delta)$ replaced with $\gamma_j^2$ respectively.
Denote $L_\epsilon$ the test  having critical region $S_\epsilon$.
\begin{lemma}\label{l1} Let $H_0$ hold.  Then the distributions of tests statistics $R^a_\epsilon(y)$ and $R_\epsilon(y)$ converge to the standard normal distribution.

For any family $\theta_\epsilon =\{\theta_{j\epsilon}\} \in Q_\epsilon$ there holds
\begin{equation}\label{i5a}
P_{\theta_\epsilon}\left(\left(T_\epsilon^a(y) - \epsilon^{-2}\rho_\epsilon - \epsilon^{-4}\sum_{j=1}^\infty \kappa_j^2\theta_{j\epsilon}^2\right)(2A_\epsilon)^{-1/2} < x_\alpha\right) = \Phi(x_\alpha)(1+o(1)).
\end{equation}
and
\begin{equation}\label{i5a}
P_{\theta_\epsilon}\left(\left(T_\epsilon(y) - C_\epsilon - \epsilon^{-4}\sum_{j=1}^\infty \kappa_j^2\theta_{j\epsilon}^2\right)(2A_\epsilon)^{-1/2} < x_\alpha\right) = \Phi(x_\alpha)(1+o(1)).
\end{equation}
as $\epsilon \to 0$.
\end{lemma}
Hence we get the following Lemma.
\begin{lemma}\label{l2} There holds
\begin{equation}\label{i6}
\beta_\epsilon(L_\epsilon) = \beta_\epsilon(L^a_\epsilon)(1+ o(1))
\end{equation}
as $\epsilon \to 0$.
\end{lemma}
\begin{lemma}\label{lx} Let $H_0$ hold.  Then the distribution of tests statistics $(T_{\epsilon\delta}(y) - C_{\epsilon\delta})(2A_\epsilon)^{-1/2}$ converge to the standard normal distribution.

There holds
\begin{equation}\label{i9a}
P_{\eta(\delta)}((T_{\epsilon\delta}(y) - C_{\epsilon\delta} - A_{\epsilon\delta})(2A_{\epsilon\delta})^{-1/2} < x_\alpha) = \Phi(x_\alpha)(1+o(1)).
\end{equation}
as $\epsilon \to 0$.
\end{lemma}
\begin{lemma}\label{l4} There holds
\begin{equation}\label{i10}
\lim_{\delta\to 0}\lim_{\epsilon \to 0} E_{\eta(\delta)} \beta_{\eta(\delta)} (L_{\epsilon\delta})  =
\lim_{\epsilon \to 0} E_{\eta_0} \beta_{\eta_0} (L_{\epsilon})
\end{equation}
where $\eta_0 = \{\eta_{0j}\}_{j=1}^ \infty$ and $\eta_{0j}$ are i.i.d. Gaussian random variables, $ E\eta_{0j} =0, \eta_{0j}^2 = \kappa_j^2, 1 \le j < \infty$.
\end{lemma}
Define Bayes a priori distribution  $P_y$ as a conditional distribution of $\eta$ given $\eta \in Q_\epsilon$. Denote $K_\epsilon = K_{\epsilon\delta}$ Bayes test with Bayes a priori distribution $P_y$. Denote $V_\epsilon$ critical region of $K_{\epsilon\delta}$.

For any sets $A$ and $B$ denote $A \triangle B = (A\setminus B)\cup(B\setminus A)$.
\begin{lemma}\label{l5} There holds
\begin{equation}\label{i11}
\lim_{\delta\to 0}\lim_{\epsilon \to 0} \int_{Q_\epsilon} P_\theta(S_{\epsilon\delta}\triangle V_{\epsilon\delta}) d\,P_y = 0
\end{equation}
and
\begin{equation}\label{i12}
\lim_{\delta\to 0}\lim_{\epsilon \to 0} P_0(S_{\epsilon\delta}\triangle V_{\epsilon\delta}) = 0.
\end{equation}
\end{lemma}
In the proof of Lemma \ref{l5} we show that the integrals in the right hand-side of (\ref{i9}) with integration domain $Q_\epsilon$ converge to one in probability as $\epsilon \to 0$. This statement is proved both for hypothesis and Bayes alternative (see \cite{er90}).

Lemmas \ref{l3}-\ref{l5} implies that, if $\alpha(K_\epsilon) = \alpha(L_\epsilon)$, then
\begin{equation}\label{i13}
\int_{Q_\epsilon} \beta_\theta(K_\epsilon) d\, P_y = \int_{Q_\epsilon} \beta_\theta(L_\epsilon) d\, P_y(1 + o(1)) = \int \beta_{\eta_0} (L_\epsilon) d\, P_{\eta_0} (1+o(1)).
\end{equation}
\begin{lemma}\label{l6} There holds
\begin{equation}\label{i14}
E_{\eta_0} \beta_{\eta_0}(L_\epsilon)  = \beta_\epsilon(L_\epsilon)(1+o(1)).
\end{equation}
\end{lemma}
Lemmas \ref{l1}, \ref{l4},  (\ref{i8}), (\ref{i13}) and Lemma \ref{l6} imply Theorem \ref{t1}.
\section{Proof of Lemmas}
Proofs of Lemmas \ref{l1},\ref{l2} and \ref{l4} are akin to the proofs of similar statements in \cite{er90} and are omitted.

\noindent{\sl Proof of Lemma \ref{l3}}. By straightforward calculations, we get
\begin{equation}\label{p1}
\sum_{j=1}^\infty E\eta_j^2(\delta) \ge \rho_\epsilon(1+\delta/2)
\end{equation}
and
\begin{equation}\label{p2}
\mbox{Var}\left(\sum_{j=1}^\infty \eta_j^2(\delta)\right) < CA_\epsilon\epsilon^{-4} \asymp  \rho_\epsilon^2 k_\epsilon^{-1}.
\end{equation}
Hence, by Chebyshev inequality, we get
\begin{equation}\label{p3}
P\left(\sum_{j=1}^\infty \eta_j^2(\delta)> \rho_\epsilon\right) = 1 +o(1)
\end{equation}
as $\epsilon \to 0$.
It remains to estimate
\begin{equation}\label{e5}
P_\mu(\eta \notin B^r_{2\infty}) = P(\max_{l_1 \le i \le l_2} i^{2r} \sum_{j = i}^{l_2} \eta_j^2-P_0(1-\delta_1/2) > P_0\delta_1/2) \le \sum_{i=l_1}^{l_2} J_i
\end{equation}
with
$$
J_i = P\left( i^{2r} \sum_{j = i}^{l_2} \eta_j^2-P_0(1-\delta_1/2)> P_0\delta_1/2\right)
$$
To estimate $J_i$ we implement the following Proposition \cite{hs}
\begin{proposition}\label{p1} Let $\xi = \{\xi_i\}_{i=1}^l$ be Gaussian random vector with i.i.d.r.v.'s $\xi_i$, $E[\xi_i] = 0, E[\xi_i^2]=1$. Let $A\in R^l\times R^l$ and $\Sigma = A^T A$. Then
\begin{equation}\label{e6}
P(||A\xi||^2 > \mbox{tr}(\Sigma) + 2\sqrt{\mbox{tr}(\Sigma^2)t} + 2 ||\Sigma||t) \le \exp\{-t\}.
\end{equation}
\end{proposition}
We put $\Sigma_i= \{\sigma_{lj}\}_{l,j=i}^{k_{\epsilon\delta}}$ with $\sigma_{jj} = j^{-2r-1}i^{2r}\frac{P_0-\delta}{2r}$ and $\sigma_{lj} =0$ if $l\ne j$.

Let $i \le k_\epsilon$. Then
\begin{equation}\label{p4}
\mbox{tr}\Sigma_i^2 = i^{4r}\sum_{j=i}^\infty \kappa_j^4(\delta) < i^{4r}((k_\epsilon -i)\kappa^4(\delta) + k_\epsilon^{-4r-1}P_0) < Ck_\epsilon^{-1}.
\end{equation}
and
\begin{equation}\label{p5}
||\Sigma_i||  \le i^{2r}\kappa^2 < Ck_\epsilon^{-1}.
\end{equation}
Therefore
\begin{equation}\label{p6}
 2\sqrt{\mbox{tr}(\Sigma_i^2)t} + 2 ||\Sigma_i||t \le C(\sqrt{k_\epsilon^{-1}t} + k_\epsilon^{-1}t)
\end{equation}
Hence, putting $t =k_\epsilon^{1/2}$, by Proposition \ref{p1}, we get
\begin{equation}\label{p7}
\sum_{i=1}^{k_\epsilon} J_i \le Ck_\epsilon \exp\{-Ck_\epsilon^{1/2}\}.
\end{equation}
Let $i \ge k_\epsilon$.  Then
\begin{equation}\label{p8}
\mbox{tr}\Sigma_i^2  < Ci^{-1}, \quad\mbox{and} \quad ||\Sigma_i|| \le Ci^{-1}
\end{equation}
Hence, putting $t =i^{1/2}$, by Proposition \ref{p1}, we get
\begin{equation}\label{p9}
\sum_{i=k_\epsilon+1}^{k_{\epsilon\delta}} J_i \le \sum_{i=k_\epsilon+1}^{k_{\epsilon\delta}} \exp\{-Ci^{1/2}\} < \exp\{-C_1k_\epsilon^{1/2}\}.
\end{equation}
Now (\ref{e5}), (\ref{p7}), (\ref{p9}) together implies Lemma \ref{l3}.

\noindent{\sl Proof of Lemma \ref{l5}}. By reasoning of the proof of Lemma 4 in \cite{er90}, Lemma \ref{l5} will be proved if we show that
\begin{equation}\label{p10}
P\left(\sum_{j=1}^\infty (\eta_j(\delta) + y_j\gamma_j(\delta)\epsilon^{-1})^2 > \rho_\epsilon\right) =1 +o(1)
\end{equation}
and
\begin{equation}\label{p11}
P\left(\sup_i i^{2r}\sum_{j=i}^\infty (\eta_j(\delta) + y_j\gamma_j(\delta)\epsilon^{-1})^2 > \rho_\epsilon\right) =1 +o(1)
\end{equation}
where $y_j, 1 \le j < \infty$  are distributed by hypothesis  or Bayes alternative.

We prove only (\ref{p11}) in the case of Bayes alternative. In other cases the reasoning are similar.

We have 
\begin{equation}\label{p12}
\begin{split}&
i^{2r}\sum_{j=i}^\infty (\eta_j(\delta) + y_j\gamma_j(\delta)\epsilon^{-1})^2=i^{2r}\sum_{j=i}^\infty \eta_j^2(\delta)\\& + i^{2r}\sum_{j=i}^\infty \eta_j(\delta)y_j\gamma_j(\delta)\epsilon^{-1} + i^{2r}\sum_{j=i}^\infty y_j^2\gamma_j^2(\delta)\epsilon^{-2} = J_{1i} + J_{2i}+ J_{3i}.
\end{split}
\end{equation}
The probability under consideration for the first addendum has been estimated in Lemma \ref{l3}.

We have
\begin{equation}\label{p13}
J_{2i} \le J_{1i}^{1/2}J_{3i}^{1/2}.
\end{equation}
Thus it remains to show that, for any $C$,
\begin{equation}\label{p14}
P_{\eta(\delta)}\left(\sup_i i^{2r}\sum_{j=i}^\infty y_j^2\gamma_j^4(\delta)\epsilon^{-2} > \delta/C\right)   =   o(1)
\end{equation}
as $\epsilon \to 0$.

Note that $y_j = \zeta_j + \epsilon\xi_j$ where $\zeta_j, y_j, 1 \le j <\infty$ are i.i.d. Gaussian random variables, $E\zeta_j =0, E\zeta_j^2 = \kappa_j^2(\delta), E \xi_j =0, E \xi_j^2 =1$.

Hence, we have
\begin{equation}\label{p15}
\begin{split}&
\sum_{j=i}^\infty y_j^2\gamma_j^4(\delta)\epsilon^{-2} =  \sum_{j=i}^\infty \gamma_j^4(\delta)\epsilon^{-2}\zeta_j^2  + \sum_{j=i}^\infty \gamma_j^4(\delta)\epsilon^{-1}\zeta_j\xi_j \\&+ \sum_{j=i}^\infty \gamma_j^4(\delta)\xi_j^2 = I_{1i}+ I_{2i}+ I_{3i}.
\end{split}
\end{equation}
Since $ \gamma_j^2\epsilon^{-2} = o(1)$, the estimates for probability of $i^{2r}I_{1i}$ are evident. It suffices to follow the estimates of (\ref{e5}). We have $I_{2i} \le I_{1i}^{1/2}I_{3i}^{1/2}$. Thus it remains to show that, for any $C$
\begin{equation}\label{p16}
P_{\eta(\delta)}\left( \sup_i i^{2r}\sum_{j=i}^\infty \gamma_j^4(\delta)\xi_j^2 > \delta/C\right) =o(1)
\end{equation}
as $\epsilon \to 0$.
Since $\gamma_j^2 = \kappa_j^2(1 + o(1)) = o(1)$, this estimate is also follows from estimates (\ref{e5}).

\noindent{\sl Proof of Lemma \ref{l6}}. By Lemmas \ref{l1}, \ref{l2} and \ref{l4}, it suffices to show that
\begin{equation}\label{p17}
\inf_\theta \sum_{j=1}^\infty \kappa_j^2\theta_j^2 =  \sum_{j=1}^\infty \kappa_j^2. 
\end{equation}
Denote $u_k = k^{2r}\sum_{j=k}^\infty \theta_j^2$. Note that $u_k \le P_0$.

Then $ \theta_j^2 = u_j j^{-2r} - u_{j+1}(j+1)^{-2r}$. Hence we have
\begin{equation}\label{p18}
\begin{split}&
A(\theta) = \sum_{j=1}^\infty \kappa_j^2\theta_j^2 = \kappa^2 \sum_{j=1}^{k_\epsilon} \theta_j^2 + \sum_{j=k_\epsilon}^\infty\kappa_j^2(u_j j^{-2r} - u_{j+1}(j+1)^{-2r})\\& =  \kappa^2 \sum_{j=1}^{k_\epsilon} \theta_j^2 + \kappa^2u_{k_\epsilon} k_\epsilon^{-2r} +  \frac{P_0}{r}\sum_{j=k_\epsilon+1}^\infty u_j (j^{-4r-1} - (j-1)^{-2r-1}j^{-2r}).
\end{split}
\end{equation}
Since $j^{-4r-1} - (j-1)^{-2r-1}j^{-2r}$ is negative, then $\inf A(\theta)$ is attained for $u_j = P_0$ and therefore $\theta_j^2 = \kappa_j^2$ for $j> k_\epsilon$.

Thus the problem is reduced  to the solution of the  following problem
\begin{equation}\label{p18}
\inf_{\theta_j} \kappa^2 \sum_{j=1}^{k_\epsilon} \theta_j^2 + \sum_{j=k_\epsilon+1}^\infty \kappa_j^4
\end{equation}
if
$$
\sum_{j=1}^{k_\epsilon} \theta_j^2 + \sum_{j=k_\epsilon+1}^\infty \kappa_j^2 =\rho_\epsilon
$$
and 
$$
k^{2r}\sum_{j=k}^\infty \theta_j^2 < P_0, \quad 1 \le j < \infty,
$$
with $\theta_j^2 = \kappa_j^2$ for $j \ge k_\epsilon$.

It is easy to see that this infimum is attained if $\kappa_j^2 = \theta_j^2 = \kappa^2$ for $j \le k_\epsilon$.

\end{document}